# A Discrete Gumbel Distribution

*Subrata Chakraborty* [a] *and Dhrubajyoti Chakravarty* [b]

[a] *Department of Statistics, Dibrugarh University, Dibrugarh, Assam, India.*

[b] *Department of Statistics, G. C. College, Silchar, Assam, India.*

email: subrata_arya@yahoo.co.in

**Abstract**

A discrete version of the Gumbel (Type I) extreme value distribution has been derived by using the general approach of discretization of a continuous distribution. Important distributional and reliability properties have been explored. It has been shown that depending on the choice of parameters the proposed distribution can be positively or negatively skewed; possess long tail(s), and exhibits equal, over or under dispersion. Log concavity of the distribution and consequential results has been established. Estimation of parameters by method of maximum likelihood, method of moments, and method of proportions has been presented. A method of checking model adequacy and regression type estimation based on empirical survival function has also been examined. Simulation study has been carried out to check the efficacy of the maximum likelihood estimators. Finally, the proposed distribution has been applied to model three real life count data regarding maximum flood discharges and annual maximum wind speeds from literature.

## 1. Introduction

### 1.1 Gumbel Distribution

The class of continuous distributions extreme value [9] has three types of extreme value distributions among which the Gumbel (Type I) distribution is the most commonly referred one, which by default is known as the extreme value distribution.

Gumbel distribution is a unimodal distribution with probability density function (pdf), cumulative distribution function (cdf), survival function (sf), failure (hazard) rate function and moment generating function (pgf) respectively given by (see [9], [23]).

$$f(x) = \sigma^{-1} e^{-(x-\mu)/\sigma} \exp[-e^{-(x-\mu)/\sigma}] \ ; -\infty < x < \infty, -\infty < \mu < \infty, \sigma > 0 \qquad (1)$$

$$F(x) = \exp[-e^{-(x-\mu)/\sigma}]$$

$$S(x) = 1 - \exp[-e^{-(x-\mu)/\sigma}]$$

$$r(x) = f(x)/S(x) = [\beta^{-1} e^{-(x-\mu)/\sigma}]/[-1 + \exp\{e^{-(x-\mu)/\sigma}\}] \text{ and}$$

$$M_X(t) = E[e^{tX}] = e^{at}\Gamma(1-bt), \ t < 1/b$$

In the case when $\mu = 0$ and $\sigma = 1$ the distribution is referred to as the standard extreme value distribution with pdf $f(x) = e^{-x} e^{-e^{-x}}$. Throughout the rest of this article a random variable $X$ following Gumbel (Type I) distribution has been referred to as $EV(\mu, \sigma)$.

In many real life situations where the variable under investigation is modeled by continuous extreme value distribution it has been observed that often the variable is actually recorded as an integer valued one instead of real valued, either because of its inherent nature or because of the limitation of measuring instruments. As such it is desirable to introduce discrete version of the existing continuous distributions.

**1.2 Discretization of Continuous Probability Distribution**

Discretization of continuous probability distribution can be done using different methodologies. One of them is as follows:

If the underlying continuous random variable $X$ has the sf $S_X(x) = \Pr(X \geq x)$ then the random variable $Y = \lfloor X \rfloor$ = largest integer less or equal to $X$ will have the probability mass function (pmf)

$$\Pr(Y = y) = \Pr(y \leq X < y+1) = \Pr(X \geq y) - \Pr(X \geq y+1)$$

$$= S_X(y) - S_X(y+1), \quad y \in Z, \text{ where } Z \text{ is the set of integers.} \tag{2}$$

So, given any continuous distribution it is possible to generate corresponding discrete distribution using the formula (2) above. The pmf of random variable $Y$ thus constructed can be viewed as discrete concentration [4] of the pdf of $X$.

Such discrete distribution retains the same functional form of the sf as that of the continuous one. As a result, many reliability characteristics remain unchanged. Discretization of many well known continuous distributions has been studied using this approach (see [21]). Notable among them are discrete Weibull distribution ([13], [15], [24], and [26]), discrete normal distribution [3], discrete Rayleigh distribution [4], discrete Maxwell distribution [7], discrete Burr and discrete Pareto distributions [8], discrete inverse Weibull distribution [12], discrete gamma distribution [20], discrete generalized gamma distribution [19] and discrete logistic distribution [21]. Besides this, there are various methods of discretization of continuous probability distributions. A partially complete list of these methods can be found in Lai [2].

The main motivation of this paper is to provide a discrete analogue of the Gumbel's extreme value distribution by discretizing the continuous Gumbel extreme value distribution in equation (1) and to investigate its important properties and applications. The proposed distribution with its structural properties has been derived in section 2. Monotonic and reliability properties have been presented in section 3 and 4 respectively. Maximum likelihood estimation (MLE) method along with method of moments and method of proportions and a simulation study for MLEs has been considered in section 5. Real life extreme value count data modeling has been presented in section 6. Summary and concluding remarks have been presented in the final section.

## 2. Discrete Gumbel Distributions

The discrete Gumbel distribution has been derived by considering the Gumbel (Type I) extreme value distribution in equation (1) using the discretization approach discussed in equation (2), after the re-parameterization $p = e^{-1/\sigma}$, and $\alpha = p^{-\mu}$.

**Definition 1**. A discrete random variable $Y$ taking values in the set of integers $Z$ is said to follow discrete Gumbel distribution with two parameters $(\alpha, p)$ that is $\text{DGUD}(\alpha, p)$, if its pmf is given by

$$f_Y(y; \alpha, p) = \Pr(Y = y) = e^{-\alpha p^{y+1}} - e^{-\alpha p^y}; y \in Z, 0 < p < 1, \alpha > 0. \qquad (3)$$

In particular,

(i) one parameter $\text{DGUD}(p)$ can be obtained for $\alpha = 1$ with pmf

$$f_Y(y; p) = \Pr(Y = y) = e^{-p^{y+1}} - e^{-p^y}; y \in Z, 0 < p < 1. \qquad (4)$$

(ii) the standard DEVD is obtained when $\alpha = 1$ and $p = \exp(-1) = 0.3679$

$$f_Y(y; p) = \Pr(Y = y) = e^{-e^{-y-1}} - e^{-e^{-y}}; y \in Z. \qquad (5)$$

## 2.1 Structural Properties of the $\text{DGUD}(\alpha, p)$

### 2.1.1 Shape of the pmf

The pmf plots of $\text{DGUD}(\alpha, p)$ for various choices of the values of the parameters have been presented in Figure 1.

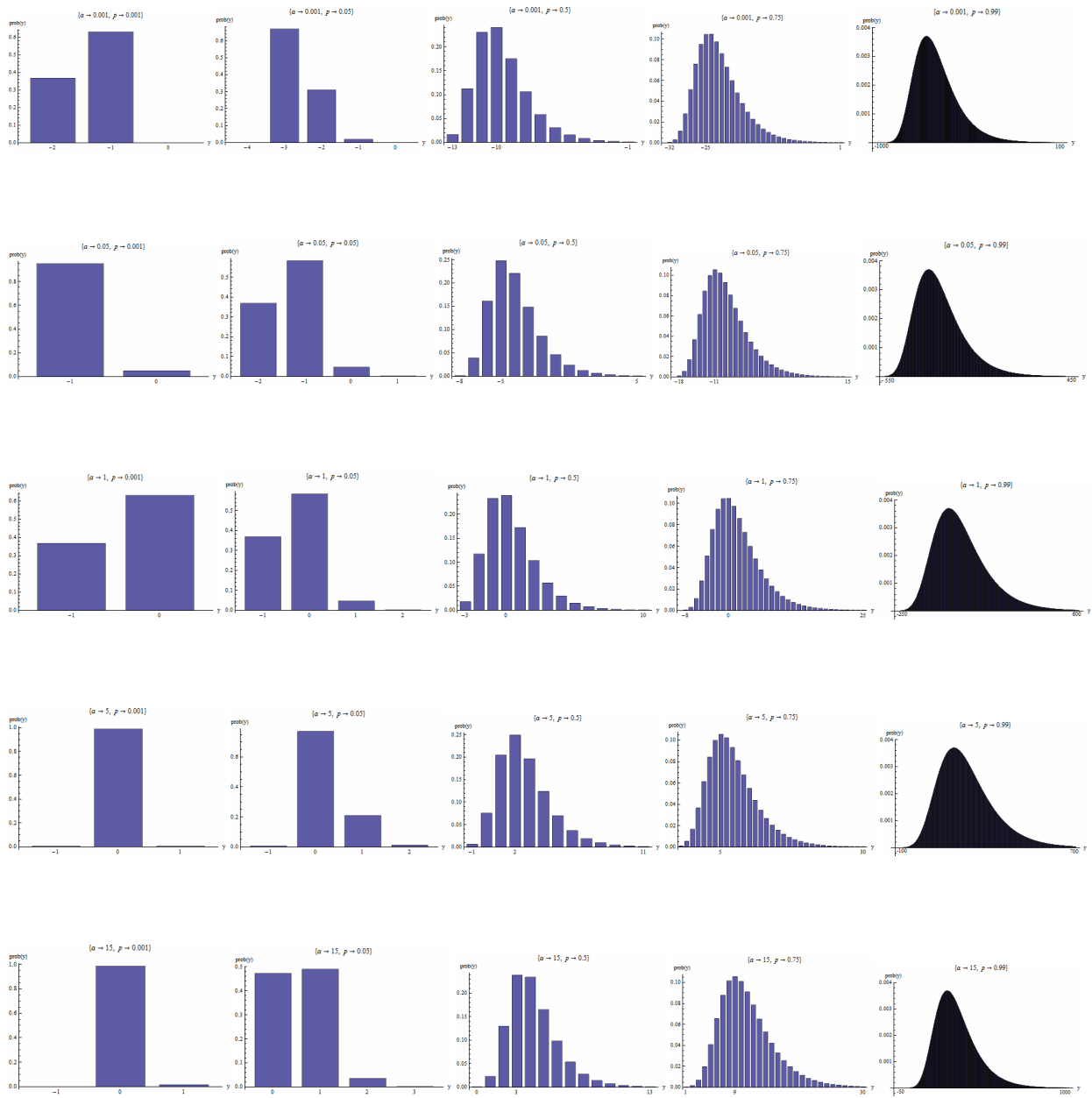

**Figure 1.** pmf of $\text{DGUD}(\alpha, p)$ for different values of $\alpha$ and $p$.

**Interpretation of the parameters of** $\text{DGUD}(\alpha, p)$: For this distribution $\alpha$ is the location and $p$ is the scale parameter. From the plots in figure 1, the following observations have been made:

(i) For moderate value of $p$, as $\alpha \to 0(\infty)$ the almost entire support tends to subset of negative (positive) integers (see section 2.1.3).

(ii) The distribution is always unimodal (see section 3.1.6 for a proof). For $\alpha=1$, the unique mode is at '0' and shifts to the right (left) from '0' as the values of the parameter $\alpha$ increases (decreases) from 1.

(iii) As $p \to 1$, the support of the distribution extends in both right and left of the mode, that is the tails are stretched with the tail to the right of the mode stretching quicker than the left (see section 2.1.2 for a proof).

### 2.1.2 Recurrence Relation for Probabilities and Long Tailedness

The probabilities can be calculated recursively using the following relation:

$$f_Y(y+1;\alpha,p) = \frac{e^{-\alpha p^{y+2}} - e^{-\alpha p^{y+1}}}{e^{-\alpha p^{y+1}} - e^{-\alpha p^{y}}} f_Y(y;\alpha,p); y \in Z, 0 < p < 1, \alpha > 0.$$

For moderate values of $\alpha$, the long tailedness of the distribution relative to Poisson distribution can be investigated by the limiting value of this ratio [17]. It is easy to see that for $y \to \infty$, the ratio of successive probabilities $f_Y(y+1;\alpha,p)/f_Y(y;\alpha,p)$ tends to $p$ as

$$\lim_{y \to \infty} \frac{e^{-\alpha p^{y+2}} - e^{-\alpha p^{y+1}}}{e^{-\alpha p^{y+1}} - e^{-\alpha p^{y}}} = \lim_{y \to \infty} \frac{pe^{-\alpha p^{y+1}} - p^2 e^{-\alpha p^{y+2}}}{e^{-\alpha p^{y}} - pe^{-\alpha p^{y+1}}} = p \text{ (using L'Hospital rule)}$$

While for $y \to -\infty$, the ratio $f_Y(y;\alpha,p)/f_Y(y+1;\alpha,p)$ of successive probabilities tends to $1/p$.

Therefore the distribution may have either right or left or both long tails depending on the value of the parameter $p$. This can be seen in the pmf plots in figure 1, where both the tails gets longer as the value of $p$ increases to unity, whereas for $p$ close to zero the right tail disappears and only the tail to the left of the mode remains.

### 2.1.3 Cumulative Distribution Function

The cdf of $DGUD(\alpha, p)$ is obtained as

$$F(y) = \Pr(Y \le y) = 1 - S(y) + \Pr(Y = y) = 1 - (1 - e^{-\alpha p^{y}}) + (e^{-\alpha p^{y+1}} - e^{-\alpha p^{y}})$$
$$= e^{-\alpha p^{y+1}}; y \in Z, 0 < p < 1, \alpha > 0. \tag{6}$$

Following results can be seen as direct consequence of equation (6):

- The proportion of positive values $= 1 - F(0) = 1 - e^{-\alpha p^{0+1}} = 1 - e^{-\alpha p}$
- The proportion of negative values $= F(-1) = e^{-\alpha p^{-1+1}} = e^{-\alpha}$
- The proportion of zeros is given by $\Pr(Y = 0) = e^{-\alpha p^{0+1}} - e^{-\alpha p^{0}} = e^{-\alpha p} - e^{-\alpha}$

- $P(a < Y \leq b) = e^{-\alpha p^{b+1}} - e^{-\alpha p^{a+1}}; a,b \in Z, 0 < p < 1, \alpha > 0.$

- Proportion of negative (non negative) values depends solely on the value of $\alpha$.

### 2.1.4 Quantiles and Random Number Generation

**Definition 2.** The point $y_u$ is known as the $u^{th}$ quantile of a discrete random variable $Y$, if it satisfies $\Pr(Y \leq y_u) \geq u$ and $\Pr(Y \geq y_u) \geq 1 - u$ that is $F(y_u - 1) < u \leq F(y_u)$. (pp. 81, Rohatgi and Saleh [25]).

**Theorem 1.** The $u^{th}$ quantile $y_u$, of $DGUD(\alpha, p)$ is given by can be obtained by $\left\lceil \dfrac{\log(1/\alpha) + \log(\log(1/u))}{\log(p)} - 1 \right\rceil$, where $\lceil a \rceil$ denotes the smallest integer greater than or equal to $a$.

***Proof.*** $\Pr(Y \leq y_u) \geq u$

$\Rightarrow e^{-\alpha p^{y_u+1}} \geq u$

$\Rightarrow -\alpha p^{y_u+1} \geq \log(u)$

$\Rightarrow p^{y_u+1} \leq -\log(u)/\alpha = \log(1/u)/\alpha$

$\Rightarrow (y_u + 1)\log(p) \leq \log((\log(1/u))/\alpha)$

$\Rightarrow y_u + 1 \geq [\log(1/\alpha) + \log(\log(1/u))] / \log(p)$ [since $o < p < 1$]

$\Rightarrow y_u \geq [[\log(1/\alpha) + \log(\log(1/u))] / \log(p)] - 1$ \hfill (7)

Similarly $\Pr(Y \geq y_u) \geq 1 - u$ gives

$y_u \leq [\log(1/\alpha) + \log(\log(1/u))] / \log(p)$ \hfill (8)

Combining (7) and (8) gives

$$\dfrac{\log(1/\alpha) + \log(\log(1/u))}{\log(p)} - 1 < y_u \leq \dfrac{\log(1/\alpha) + \log(\log(1/u))}{\log(p)}$$

Hence $y_u$ is an integer given by $y_u = \left\lceil \dfrac{\log(1/\alpha) + \log(\log(1/u))}{\log(p)} - 1 \right\rceil$. \hfill (9)

The three quartiles of $DEVD(\alpha, p)$ can be computed easily by the following formulae:

$$y_{0.025} = \left\lceil \dfrac{\log(1/\alpha) + 0.3266}{\log(p)} - 1 \right\rceil$$

$$y_{0.50} = \left\lceil \frac{\log(1/\alpha) - 0.3665}{\log(p)} - 1 \right\rceil$$

$$y_{0.75} = \left\lceil \frac{\log(1/\alpha) - 1.2470}{\log(p)} - 1 \right\rceil$$

A random number (integer) can be sampled from the proposed model through the usual inverse transformation method. Let $U$ be a random number drawn from a uniform distribution on (0, 1); then a random number $Y$ following $\text{DGUD}(\alpha, p)$ is obtained computing the right hand side of equation (9).

### 2.1.5 Mode

**Theorem 2.** Mode at $y = \lfloor -\log(\alpha)/\log(p) \rfloor$

*Proof.* Since $\alpha = p^{-\mu}$ $(\alpha > 0)$, the pmf of $\text{DGUD}(\alpha, p)$ can be rewritten as

$$f_Y(y; \alpha, p) = \Pr(Y = y) = e^{-\alpha p^{y+1}} - e^{-\alpha p^y} = e^{-p^{y-\mu+1}} - e^{-p^{y-\mu}}, \quad -\infty < \mu < \infty.$$

Now, for any value of $\mu$,

$$\Pr(Y = \lfloor \mu \rfloor) - \Pr(Y = \lfloor \mu \rfloor + 1) = e^{-p^1} - e^{-p^0} - e^{-p^2} + e^{-p^1} = 2e^{-p} - (e^{-1} - e^{-p^2}),$$

which is always positive for $0 < p < 1$.

Similarly,

$$\Pr(Y = \lfloor \mu \rfloor) - \Pr(Y = \lfloor \mu \rfloor - 1) = e^{-p^1} - e^{-p^0} - e^{-p^0} + e^{-p^{-1}} = e^{-p} - 2/e + e^{-(1/p)},$$

which is also positive for $0 < p < 1$.

But as will be been seen later in corollary 4 in section 3 that the distribution is unimodal hence the only mode is at $Y = \lfloor \mu \rfloor = \lfloor -\log(\alpha)/\log(p) \rfloor$.

This can be easily verified with pmf plots given in figure 1.

### 2.1.6 Moments

The $r^{th}$ moment about origin for is given by

$$E(Y^r) = \sum_{y=-\infty}^{\infty} y^r (e^{-\alpha p^{y+1}} - e^{-\alpha p^y})$$

There is no close form for the moments. The values of mean and variance for different combination of parameter values have been plotted in Figure 2 to assess presence of any pattern.

Contours of mean and variance are plotted in Figure 3 and 4 respectively to partition the parameter space in different regions according to the values of mean and variance respectively.

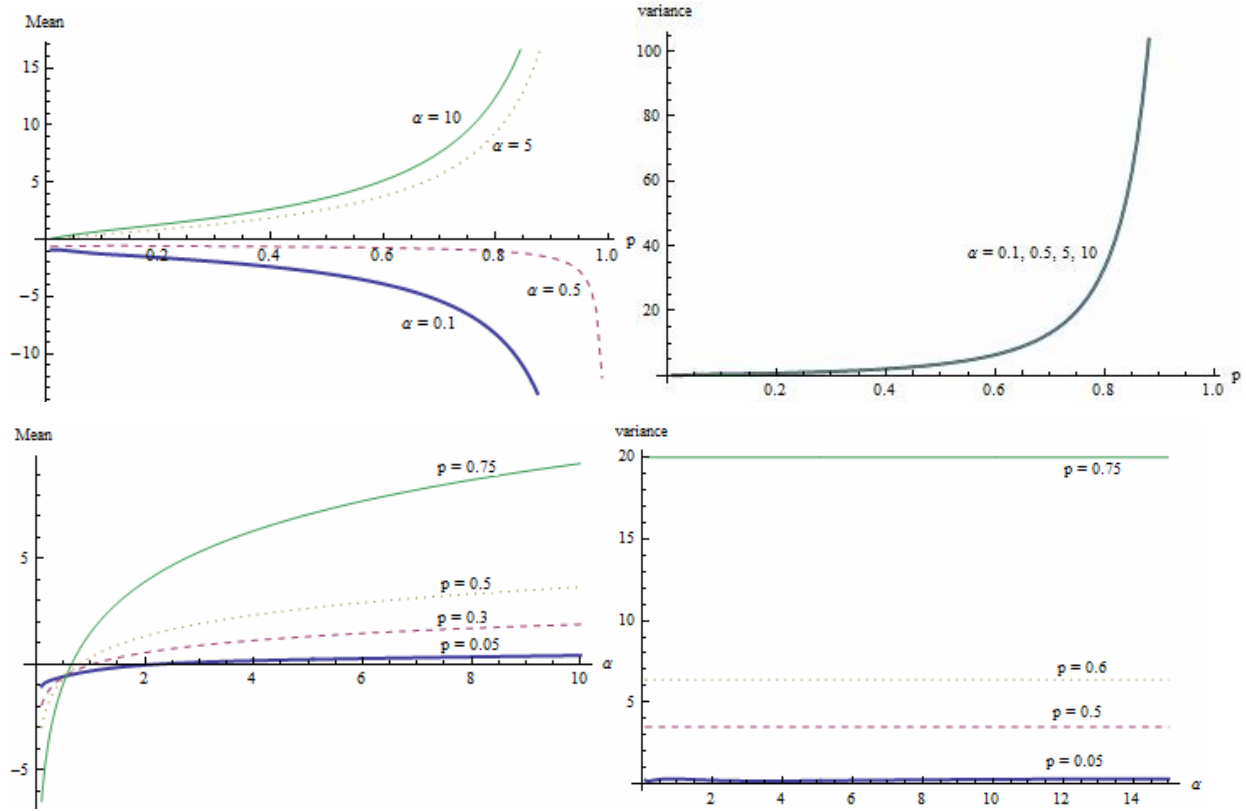

**Figure 2.** Mean and Variance of $\text{DGUD}(\alpha, p)$.

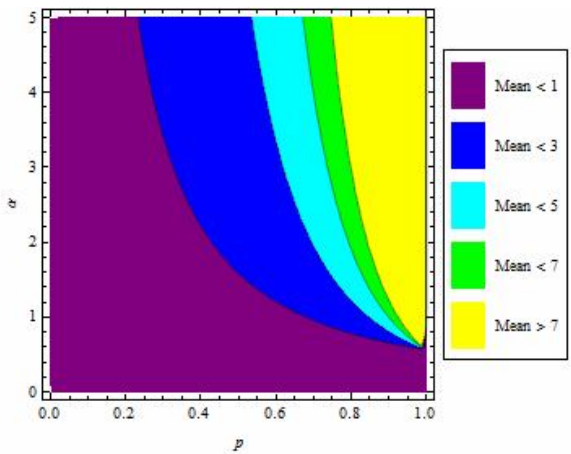

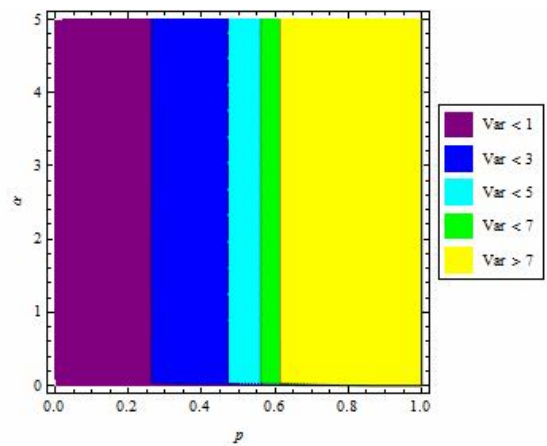

**Figure 3.** Mean Contour plot

**Figure 4.** Variance Contour plot

From the plots in figures 2-4 it is apparent that both mean and the variance of $\text{DGUD}(\alpha, p)$ have bounds and that the variance is affected only by the value of parameter $p$ not by the value of $\alpha$. These findings have also been established mathematically in the theorem 3.

**Theorem 3**. For $DGUD(\alpha, p)$ mean and variance are both bounded as

$\log \alpha / \log(1/p) + 0.577216/\log(1/p) - 1 \leq E(Y) \leq \log \alpha / \log(1/p) + 0.577216/\log(1/p)$ and

$(\pi^2/6)[1/\log(1/p)] \leq Var(Y) \leq (\pi^2/6)[1/\log(1/p)]^2 + 1/4$

***Proof***. For continuous random variable $X$ following $EV(\mu, \sigma)$ distribution with pdf in (1) it is known that (see [9], [14], and [23])

$E(X) = \mu + 0.577216\,\sigma = \log\alpha/\log(1/p) + 0.577216/\log(1/p)$ and

$Var(X) = \pi^2\sigma^2/6 = (\pi^2/6)[1/\log(1/p)]^2$, where $p = \exp(-1/\sigma)$ and $\alpha = p^{-\mu}$.

Now, the discretized version $Y$ of $X$, that is, $DGUD(\alpha, p)$ is defined as $Y = \lfloor X \rfloor$ = largest integer less or equal to $X$. It can be assumed that $X = Y + U$, where $U = X - \lfloor X \rfloor$, is the fractional part of $X$ which is chopped off from $X$ to obtain $Y$. Hence $U$ will also have a pdf in the support $(0, 1)$ and independent of $Y$. Therefore

$E(Y) = E(X - U) = E(X) - E(U)$

$= \log\alpha/\log(1/p) + 0.577216/\log(1/p) - E(U)$

But for any continuous random variable $U$ in $(0,1)$, $0 < E(U) \leq 1$, hence the result follows. By similar argument,

$Var(Y) = Var(X - U) = Var(X) + Var(U)$ [Assuming independence of $X$ and $U$]

$= (\pi^2/6)[1/\log(1/p)]^2 + Var(U)$.

But for any continuous random variable $U$ in $(0,1)$, $0 < Var(U) \leq 1/4$, hence the result follows.

**Corollary 1**: It is obvious that for $\alpha > (<) 1$ the mean is an increasing (decreasing) function of $p$, while the variance is independent of $\alpha$ and is an increasing function of $p$.

**Corollary 2**: $E(Y)$ and $Var(Y)$ can be respectively approximated by the formulae

$E(Y) \cong \log\alpha/\log(1/p) + 0.577216/\log(1/p) - 0.5$

and $Var(Y) \cong (\pi^2/6)[1/\log(1/p)]^2 + 0.125$

### 2.1.7 Skewness and Kurtosis

There are a number of measures of skewness in the literature to quantify the lack of symmetry in a probability distribution among them the most commonly used one are Pearson's measure defined by ((mean - mode) / standard deviation), Yule's measure of skewness $\gamma_1$ which is defined as the third central moment divided by second central moment raised to the power 3/2.

But there is no single measure which has the universal acceptability [22]. So instead of using the common measure here a new measure of homogeneous skeweness proposed by Das et al.[22] is considered for investigating the skeweness of $\mathrm{DGUD}(\alpha, p)$.

**Definition 3.** The measure of homogeneous skeweness (HSK) for a discrete unimodal distribution with pmf $f_Y(y)$ and mode at *M* is defined as [22]

$$\mathrm{HSK} = \sum_{y>0}\{f_Y(M+y) - f_Y(M-y)\} + f_Y(M)*\mathrm{Sign}(\sum_{y>0}\{f_Y(M+y) - f_Y(M-y)\})$$

*Interpretation*:

HSK=0 if and only if the distribution is symmetric.

HSK=1(-1) if and only if the distribution is decreasing (increasing).

HSK $>$ ($<$) 0 if and only if the distribution is homogeneously right (left) skewed.

**Coefficient of kurtosis** $\beta_2$ is defined as the fourth central divided by the square of the second central moment of the probability distribution. These coefficients are vital for knowing the shape of a distribution.

For our model, large and complicated expressions have been obtained for these indices. Table 1 shows some values of these quantities for different values of the parameters $\alpha$ and $p$. From the table 1 (also in figure 1), it can be seen that $\mathrm{DGUD}(\alpha, p)$ may be positively skewed (right tail is longer than the left tail) or negatively skewed (right tail is shorter than the left tail) depending upon the choice of the values of the parameters, which indicates the potentiality of application of this distribution to differently skewed real life data. The distribution has been seen as leptokurtic in most cases except for the cases when either $\alpha$ or $p$ or both are very close to zero.

**Table 1. Skewness and Kurtosis (in parenthesis) of the $\mathrm{DGUD}(\alpha, p)$ for various choice of $\alpha$ and $p$.**

| $\alpha \rightarrow$ $p \downarrow$ | 0.001 | 0.050 | 0.150 | 0.500 | 5.00 | 15.0 |
|---|---|---|---|---|---|---|
| 0.001 | -.998(1.349) | 1.000(18.75) | 1.000(5.420) | 1.000(1.235) | -0.990(85.57) | 1.000(65.80) |
| 0.100 | -.810(3.793) | 0.987(4.364) | -.721(4.389) | 0.987(4.364) | 0.987(4.364) | -.721(4.389) |
| 0.250 | -.548(4.971) | 0.919(5.026) | 0.819(5.028) | 0.729(4.990) | -0.463(4.935) | 0.953(4.993) |
| 0.400 | 0.565(5.207) | 0.716(5.205) | 0.808(5.205) | 0.427(5.207) | 0.729(5.205) | 0.819(5.205) |
| 0.500 | 0.282(5.287) | 0.596(5.287) | 0.398(5.287) | 0.264(5.287) | 0.427(5.287) | 0.693(5.287) |
| 0.650 | 0.561(5.355) | 0.279(5.355) | 0.451(5.355) | 0.388(5.355) | 0.493(5.355) | 0.355(5.355) |
| 0.750 | 0.470(5.380) | 0.388(5.380) | 0.350(5.380) | 0.389(5.380) | 0.389(5.380) | 0.352(5.380) |
| 0.850 | 0.323(5.394) | 0.332(5.394) | 0.303(5.394) | 0.352(5.394) | 0.372(5.394) | 0.343(5.394) |
| 0.900 | 0.298(5.397) | 0.308(5.397) | 0.341(5.397) | 0.297(5.397) | 0.286(5.397) | 0.319(5.397) |
| 0.990 | 0.269(6.308) | 0.271(4.855) | 0.266(4.502) | 0.264(3.973) | 0.265(2.695) | 0.268(3.057) |

From the table 1 it can be seen that the distribution may be decreasing or increasing or right skewed or left skewed. (See also the figure 1 for a visual verification these findings for example for $\alpha = 0.001, p = 0.001$ in the plot in first row of first column the pmf is increasing since HSK is almost -1, when $\alpha = 0.05, p = 0.001$, in the plot in second row of first column the pmf is decreasing since HSK is almost +1 and, when $\alpha = 15, p = 0.5$, in the 3rd plot in last row the pmf is right skewed since HSK is 0.693 etc.).

## 3. Monotonic Properties

Log concavity is an important property of probability distribution. Characteristic like reliability function, distribution function, failure rate function, mean residual life time function, tail probabilities and moments of a log concave probability distribution have specific properties (see [1, 35]). Examples of discrete log concave probability distributions include Bernoulli, binomial, Poisson, geometric and negative binomial etc. whereas the logarithmic series distribution is a non concave discrete distribution. Here the log concavity of $DGUD(\alpha, p)$ distribution along with implications has been presented.

**Definition 3**. A distribution with pmf $f_Y(y)$ is log concave ([1], [6], [16] and [29]) if

$$f_Y(y+1)^2 > f_Y(y) f_Y(y+2).$$

**Definition 4**. A discrete probability distribution $p_k = P(X = k)$ with support on the lattice of integers is unimodal [10] if there exists at least one integer $M$ such that

$$p_k \geq p_{k-1} \; \forall \, k \leq M \text{ and } p_{k+1} \geq p_k \; \forall \, k \geq M.$$

A unimodal distribution is called strongly unimodal if its convolution with any other unimodal distribution is again unimodal. Clearly strong unimodality implies unimodality. Some examples of strongly unimodal distributions are the binomial and Poisson [10]. Strongly unimodal distribution therefore may have one (unique) or more modes.

**Theorem 4.** The $DGUD(\alpha, p)$ is log-concave.

*Proof.* To prove the log–concavity of the distribution it suffices to prove that

$$\Rightarrow \{e^{-\alpha p^{y+1}} - e^{-\alpha p^y}\}^2 \geq \{e^{-\alpha p^{y+2}} - e^{-\alpha p^{y+1}}\}\{e^{-\alpha p^y} - e^{-\alpha p^{y-1}}\}, 0 < p < 1, \alpha > 0, y \in Z.$$

Let $g(t) = e^{-\alpha p^t}$, where $\alpha > 0, 0 < p < 1$. Then, $0 \leq g(t) \leq 1$

Since $g'(t) = \alpha \log(1/p) \, p^t \, e^{-\alpha p^t} > 0 \; \forall \, t \in Z$, it is an increasing function of $t$.

The root of the equation $g''(t) = 0$ is $\log(1/\alpha)/\log(p)$. Let $T = \log(1/\alpha)/\log(p)$.

Also $g''(t) > 0\ \forall t \in (-\infty, T)$, implying that $g(t)$ is convex for $t \in (-\infty, T)$. (10)

$g''(t) < 0\ \forall t \in (T, \infty)$, implying that $g(t)$ is concave for $t \in (T, \infty)$ (11)

Hence, for every $t$, $s \leq T$, $g((t+s)) \leq (g(t) + g(s))/2$ and for every $t$, $s \geq T$,

$g((t+s)) \geq (g(t) + g(s))/2$. Also, it implies that for every $t < s$, $0 \leq g(s) - g(t) \leq 1$.

For simplification in notations, let

$A_0 = g(y), A_1 = g(y+1), A_2 = g(y+2), A_{(-1)} = g(y-1)$. Obviously,

$A_{(-1)} \leq A_0 \leq A_1 \leq, A_2$.

With these notations, it is to be proved that $(A_1 - A_0)^2 \geq (A_2 - A_1)(A_0 - A_{-1})$.

Since $0 \leq (A_2 - A_1) \leq 1$ and $0 \leq (A_0 - A_{(-1)}) \leq 1$, it suffices to prove one of the following:

$A_1 - A_0 \geq A_2 - A_1$ and $A_1 - A_0 \geq A_0 - A_{(-1)}$. This can be proved in three cases as below:

Let us take two values as $t = y-1$ and $s = y+1$.

**Case I.** If $y+1 \leq T$, then from (10), it is seen that

$2A_0 \leq A_1 - A_{(-1)} \Rightarrow A_1 - A_0 \geq A_0 - A_{(-1)}$.

**Case II.** If $y \geq T$, then from (11), it is seen that

$2A_0 \leq A_1 - A_{(-1)} \Rightarrow A_1 - A_0 \geq A_0 - A_{(-1)}$.

**Case III.** For $y \leq T \leq y+1$, by convexity and concavity of $g(t) = e^{-\alpha p^t}$ in the intervals $(y-1, T)$ and $(y-1, T)$ respectively, one can get

$A_1 \geq kA_2 + (1-k)g(T); 0 < k < 1$ and $A_0 \leq kA_{(-1)} + (1-k)g(T); 0 < k < 1$

Simplifying the inequalities of above three cases leads to the desired result.

**Corollary 3.** As a direct consequence of log–concavity (proposition 10 of Mark [29], also see [1, 14]), the following result holds for the proposed $DGUD(\alpha, p)$. It

➤ strongly unimodal
➤ Has an increasing failure (hazard) rate distribution.
➤ Monotonically decreasing Mean residual life function.
➤ Have all its moments.
➤ Remains log concave if truncated.

- Gives unimodal and log-concave distribution when convoluted with any other discrete distribution.

## 4. Reliability Properties

- **Survival Function (sf)**

$S(y) = \Pr(Y \geq y)$ of $\mathrm{DGUD}(\alpha, p)$ is given by

$$S(y) = \begin{cases} 1 - e^{-\alpha p^y}, & y \in Z \\ 1 - e^{-\alpha p^{\lfloor y \rfloor + 1}}, & y \in R \end{cases}.$$

- **Failure (Hazard) Rate Function**

$$r(y) = \frac{\Pr(Y = y)}{S(y)} = \frac{1 - e^{\alpha(1-p)p^y}}{1 - e^{\alpha p^y}}$$

The failure rate for any value of $\alpha (> 0)$ and $0 < p < 1$ is increasing. The result also follows from log-concavity (see theorem 4).

- **Second Failure Rate (SRF) Function**

The second failure rate defined by $r^*(y) = \log[S(y)/S(y+1)]$ was introduced by Roy and Gupta [5] (also discussed by Xie et al. [11]). For $Y \sim \mathrm{DGUD}(\alpha, p)$ the SRF is given by

$$r^*(y) = \log[(1 - e^{-\alpha p^y})/(1 - e^{-\alpha p^{y+1}})]$$

**Theorem 7.** If $X_i$'s $(i = 1, 2, 3, \cdots, n)$ be independent integer valued random variables and $Y = \mathrm{Max}(X_1, X_2, \cdots, X_n)$, then $Y$ is $\mathrm{DGUD}(\alpha, p)$, if $X_i$'s are $\mathrm{DGUD}(\alpha_i, p)$, where

$$\alpha = \sum_{i=1}^{n} \alpha_i.$$

*Proof.* Let $X_i$ $(i = 1, 2, 3 \ldots n)$ be independently distributed $\mathrm{DGUD}(\alpha_i, p)$, then

$$F(x) = e^{-\alpha_i p^{x+1}}; \ x \in Z$$

Consider, for all $y \in Z$.

$$F(y) = \Pr[Y \leq y] = \prod_{i=1}^{n} \Pr[X_i \leq y] = \prod_{i=1}^{n} e^{-\alpha_i p^{y+1}}$$

$$= (e^{-\alpha_1 p^{y+1}})(e^{-\alpha_2 p^{y+1}}) \cdots (e^{-\alpha_n p^{y+1}})$$

$$= e^{-(\alpha_1 + \alpha_2 + \cdots \alpha_n) p^{y+1}} = e^{-\alpha p^{y+1}}. \text{ Thus } Y \sim \mathrm{DGUD}(\alpha, p).$$

**Corollary 4.** If $X_i$'s, $(i = 1, 2, 3....n)$ are identically independently distributed $\text{DGUD}(\alpha, p)$ random variable and $Y = \text{Max}(X_1, X_2, \cdots, X_n)$, then $Y$ is $\text{DGUD}(n\alpha, p)$.

**Remark 1.** This property is useful as it will allow modeling reliability of a parallel system with identical components having $\text{DGUD}(\alpha, p)$.

**Theorem 6.** If $X \sim \text{EV}(\mu, \sigma)$ distribution, then $Y = \lfloor X \rfloor \sim \text{DGUD}(\alpha, p)$ with $p = e^{-1/\sigma}$ and $\alpha = p^{-\mu}$.

*Proof.* Consider,

$$\Pr[Y \geq y] = \Pr[\lfloor X \rfloor \geq y] = \Pr[X \geq y] = 1 - \exp[-e^{-(y-\mu)/\sigma}] = 1 - e^{-\alpha p^y}$$

Thus $Y = \lfloor X \rfloor \sim \text{DGUD}(\alpha, p)$.

**Theorem 7.** If $X \sim \text{EV}(\mu, \sigma)$ distribution, then $Y = \lfloor e^{-(x-\mu)/\sigma} \rfloor \sim \text{Geometric}(q)$ distribution with $q = \exp[-e^{-1/\sigma}]$, $(0 < q < 1/e)$.

*Proof.* Under the transformation $Y = \lfloor e^{-(x-\mu)/\sigma} \rfloor$, $Y$ can take values $0, 1, 2 \ldots$

Consider, $y = 0, 1, 2, \cdots$. Now,

$$\Pr[Y \geq y] = \Pr[\lfloor e^{-(X-\mu)/\sigma} \rfloor \geq y] = \Pr[e^{-(X-\mu)/\sigma} \geq y]$$

$$= \Pr[-(X-\mu)/\sigma \geq \log y] = \Pr[X \leq \mu - \sigma \log y] = \exp[e^{-(\mu - \sigma \log y + 1 - \mu)/\sigma}]$$

$$= \exp[v] = e^{-py} = (e^{-p})^y = q^y$$

Which is the survival function of a geometric r.v. Hence the result.

**Theorem 8.** If $X \sim \text{Geometric}(e^{-p})$, then $Y = \mu - \sigma \log X \sim \text{DGUD}(\alpha, p)$ with $p = e^{-1/\sigma}$ and $\alpha = p^{-\mu}$.

*Proof.* Consider, $y = 0, 1, 2, \cdots$. Now,

$$\Pr[Y \geq y] = \Pr[\mu - \sigma \log X \geq y] = \Pr[-\sigma \log X \geq y - \mu]$$

$$= \Pr[X \leq e^{(y-\mu)/\sigma}] = 1 - (e^{-p})^{e^{(y-\mu)/\sigma}} = 1 - e^{-e^{-1/\sigma} e^{(y-\mu)/\sigma}}$$

$$= 1 - e^{-e^{-(y-\mu)/\sigma}} = 1 - e^{-(e^{-1/\sigma})^y (e^{-1/\sigma})^{-\mu}} = 1 - e^{-\alpha p^y},$$

which the survival function of $\text{DGUD}(\alpha, p)$. Hence proved.

## 5. Estimation of Parameters and Simulation

In this section, the three commonly used methods for estimating the parameters have been discussed. Further a method based on empirical survival function for checking model adequacy and parameter estimation has also been discussed.

## 5.1 Maximum Likelihood Estimation (MLE)

Let $\mathbf{y} = (y_1, y_2, \ldots, y_n)$ is a random sample of size $n$ drawn from the $\mathrm{DGUD}(\alpha, p)$ distribution (3). The log-likelihood function is given by

$$\log L = \sum_{i=1}^{n} \log(e^{-\alpha p^{y_i+1}} - e^{-\alpha p^{y_i}})$$

Since there is no close form solution of the log-likelihood equations, the MLEs may be obtained by numerical method or direct numerical search for the global maximum of the log-likelihood surface.

The asymptotic variance-covariance matrix of the MLEs of parameters $\alpha$ and $p$ may be obtained from inverse of the Fishers information matrix

## 5.2 Method of Moments

Moment estimates of $\alpha$ and $p$ can be obtained by either numerically solving the equations $\sum_{y=-\infty}^{\infty} y(e^{-\alpha p^{y+1}} - e^{-\alpha p^{y}}) = m_1$ and $\sum_{y=-\infty}^{\infty} y^2(e^{-\alpha p^{y+1}} - e^{-\alpha p^{y}}) = m_2$, or by minimizing

$$\left\{\sum_{y=-\infty}^{\infty} y(e^{-\alpha p^{y+1}} - e^{-\alpha p^{y}}) - m_1\right\}^2 + \left\{\sum_{y=-\infty}^{\infty} y^2(e^{-\alpha p^{y+1}} - e^{-\alpha p^{y}}) - m_2\right\}^2$$

with respect to $\alpha$ and $p$, where $m_1$ and $m_2$ are respectively the 1st and 2nd observed sample raw moments.

Since sample moments are consistent estimators of the population moments by invariance property it can be seen that the moment estimators of $\alpha$ and $p$ obtained above are also consistent.

## 5.3 Method of Proportions

Suppose that the observed proportion of zeros, positive values and negative values in the sample be $p_0$, $p_+$ and $p_-$ respectively. Obviously, $p_0 + p_- + p_+ = 1$.

Here the parameters $\alpha$ and $p$ are estimated by simultaneously solving the two equations obtained by equating $p_+$ and $p_-$ with corresponding theoretical values given section 2.1.3. That is by solving $1 - e^{-\alpha p} = p_+$ and $e^{-\alpha p} - e^{-\alpha} = p_-$, which yields $\hat{\alpha} = \log(1/p_-)$ and $\hat{p} = \log(p_-)\log(1-p_+)$, provided that $p_- \neq 0$.

This method yields analytical estimators for the parameters unlike the MLE and MM. Also since $p_+$ and $p_-$ are unbiased and consistent estimators of corresponding population proportions therefore, the estimators of $\alpha$ and $p$ are also consistent.

## 5.4 Estimation of $p$ from Empirical Survival Function

For large samples, the empirical estimator $\hat{S}(y)$ of $S(y)$ is defined as the proportion of observations greater than or equal to $y$. In case of $\text{DGUD}(\alpha, p)$,

$$\prod_{i=1}^{n} S(y_i) = (1 - e^{-\alpha p^{y_i}})^n$$

So that when $\alpha$ is known, an empirical estimator of $p$ can be obtained as

$$\hat{p} = \exp[\{\log \alpha + \log[\log[1 - \{\prod \hat{s}(y_i)\}^{1/n}]]/y_i]$$

Since the empirical survival function $\hat{S}(y)$ is a consistent estimator of survival function $S(y)$, therefore by invariance property of consistent estimator the above estimator of p is also consistent.

## 5.5 Model diagnostic

Consider the survival function of $\text{DGUD}(\alpha, p)$:

$$S(y) = 1 - e^{-\alpha p^y} \Rightarrow 1 - S(y) = e^{-\alpha p^y} \Rightarrow -\log(1 - S(y)) = \alpha p^y$$
$$\Rightarrow -\log(1 - S(y)) = \alpha p^y \Rightarrow -\log(-\log(1 - S(y))) = -\log(\alpha) - \log(p) y$$
$$\Rightarrow z = a + by, \tag{12}$$

where $z = -\log(-\log(1 - S(y)))$, $a = -\log(\alpha) = \log(1/\alpha)$ and $b = -\log(p) = \log(1/p)$.

The linear equation in (12) serves as an important tool for checking model adequacy. By computing the empirical survival function $\hat{S}(y)$ from the data and plotting $-\log(-\log(1 - \hat{S}(y)))$ against $y$, one can prescribe $\text{DGUD}(\alpha, p)$ as an adequate model for the given data, provided the plot is nearly a straight line.

**Remark 2**. The parameters $\alpha$ and $p$ can be estimated by considering (12) as a linear regression equation and first estimating the constants $a$ and $b$ by using ordinary least square method and then by $\hat{\alpha} = e^{-a}$ and $\hat{p} = e^{-b}$.

### 5.6 Simulation Study

A thorough simulation analysis is carried out by generating $n = 1000$ replications each of samples of different sizes $k = 25, 50, 100$ for each pair $(\alpha, p)$ of values (0.05, 0.25), (0.05, 0.5), (0.05, 0.75); (1, 0.25), (1, 0.5), (1, 0.75) and (5, 0.25), (5, 0.5), (5, 0.75). MLEs of $\hat{\alpha}$ and $\hat{p}$ are obtained by applying global numerical optimization method along with the observed Fisher's information matrix for each of these samples. It may be noted that the sample is obtained by generating continuous extreme value distribution and then taking floor. The results of the simulation analysis are given in Table 2. In the table, all the entries are the means of estimates for 1000 samples. Accuracy and precision of the estimation method are checked and established using the following criteria:

I. The expected value of the estimator: $E(\hat{\theta}) = (1/n)\sum_{i=1}^{n} \hat{\theta}_i = \bar{\theta}$,

II. Average of the biases of estimates: $\text{Bias} = (1/n)\sum_{i=1}^{n}(\hat{\theta}_i - \theta)$,

III. Estimate of the standard error: $E[SE(\hat{\theta})] = \sqrt{\frac{1}{n}\sum_{i=1}^{n}\left(\frac{-\partial^2 \log L}{\partial \hat{\theta}_i^2}\right)}$,

IV. Average width of the 95% confidence interval estimators:

$$\text{AW}(\theta) = \frac{1}{n}\sum_{i=1}^{n}\left[\left(\hat{\theta}_i + 1.96\sqrt{\left(\frac{-\partial^2 \log L}{\partial \hat{\theta}_i^2}\right)}\right) - \left(\hat{\theta}_i - 1.96\sqrt{\left(\frac{-\partial^2 \log L}{\partial \hat{\theta}_i^2}\right)}\right)\right],$$

V. Coverage rate: $\text{CR}(\theta) = \text{Probability of } \theta_i \in \left(\hat{\theta}_i \pm 1.96\sqrt{\left(\frac{-\partial^2 \log L}{\partial \hat{\theta}_i^2}\right)}\right)$

where $\hat{\theta}_i$ is the estimate of the unknown true value $\theta$ obtained from the $i^{\text{th}}$ sample $(i = 1, 2, \cdots, n)$. From table 2 it can be seen that, as sample size increases the estimates of bias, standard error decreases. As the value of $\alpha$ increases, the bias in $\hat{\alpha}$ increases, whereas the bias in $\hat{p}$ decreases with the increasing value of $p$. As expected, the average width of the interval estimators decreases as sample size increases. It is seen that for class interval based on fisher's information matrix, the coverage rate that is never smaller than 0.861, and occasionally even higher than the nominal value 0.95. Except in a few cases the coverage rate increases with the sample size.

**Table 2: Simulation results: Bias and SE of point estimators, AW and CR of confidence intervals**

| $\alpha$ | $p$ | $n$ | $E(\hat{\alpha})$ | Bias$(\hat{\alpha})$ | $E[SE(\hat{\alpha})]$ | AW$(\alpha)$ | CR$(\alpha)$ | $E(\hat{p})$ | Bias$(\hat{p})$ | $E[SE(\hat{p})]$ | AW$(p)$ | CR$(p)$ | $E[Cov(\hat{\alpha},\hat{p})]$ |
|---|---|---|---|---|---|---|---|---|---|---|---|---|---|
| 0.05 | 0.25 | 25 | 0.050 | 0.000 | 0.032 | 0.116 | 0.861 | 0.234 | 0.016 | 0.060 | 0.236 | 0.908 | 0.002 |
| | | 50 | 0.050 | 0.000 | 0.022 | 0.083 | 0.899 | 0.242 | 0.008 | 0.043 | 0.167 | 0.935 | 0.001 |
| | | 100 | 0.050 | 0.000 | 0.016 | 0.059 | 0.917 | 0.247 | 0.003 | 0.030 | 0.118 | 0.941 | 0.000 |
| 0.05 | 0.50 | 25 | 0.050 | 0.000 | 0.031 | 0.110 | 0.866 | 0.482 | 0.018 | 0.057 | 0.225 | 0.937 | 0.002 |
| | | 50 | 0.050 | 0.000 | 0.021 | 0.079 | 0.900 | 0.500 | 0.000 | 0.040 | 0.158 | 0.949 | 0.001 |
| | | 100 | 0.050 | 0.000 | 0.015 | 0.057 | 0.920 | 0.495 | 0.005 | 0.028 | 0.111 | 0.943 | 0.000 |
| 0.05 | 0.75 | 25 | 0.050 | 0.000 | 0.030 | 0.109 | 0.869 | 0.738 | 0.012 | 0.036 | 0.139 | 0.940 | 0.001 |
| | | 50 | 0.049 | 0.001 | 0.021 | 0.077 | 0.909 | 0.743 | 0.007 | 0.025 | 0.096 | 0.957 | 0.000 |
| | | 100 | 0.050 | 0.000 | 0.015 | 0.055 | 0.919 | 0.747 | 0.003 | 0.017 | 0.067 | 0.943 | 0.000 |
| 1 | 0.25 | 25 | 1.056 | -0.056 | 0.257 | 0.978 | 0.948 | 0.231 | 0.019 | 0.062 | 0.240 | 0.898 | 0.001 |
| | | 50 | 1.019 | -0.019 | 0.170 | 0.659 | 0.938 | 0.240 | 0.010 | 0.043 | 0.170 | 0.941 | 0.001 |
| | | 100 | 1.009 | -0.009 | 0.117 | 0.457 | 0.940 | 0.245 | 0.005 | 0.031 | 0.120 | 0.938 | 0.001 |
| 1 | 0.50 | 25 | 1.035 | -0.035 | 0.229 | 0.879 | 0.943 | 0.483 | 0.018 | 0.057 | 0.225 | 0.931 | 0.003 |
| | | 50 | 1.012 | -0.012 | 0.156 | 0.606 | 0.942 | 0.490 | 0.010 | 0.040 | 0.158 | 0.950 | 0.002 |
| | | 100 | 1.004 | -0.004 | 0.109 | 0.425 | 0.936 | 0.495 | 0.005 | 0.028 | 0.111 | 0.940 | 0.001 |
| 1 | 0.75 | 25 | 1.035 | -0.035 | 0.223 | 0.859 | 0.934 | 0.737 | 0.013 | 0.036 | 0.139 | 0.936 | 0.002 |
| | | 50 | 1.011 | -0.011 | 0.152 | 0.593 | 0.941 | 0.743 | 0.007 | 0.025 | 0.097 | 0.953 | 0.001 |
| | | 100 | 1.004 | -0.004 | 0.107 | 0.416 | 0.934 | 0.747 | 0.003 | 0.017 | 0.067 | 0.950 | 0.001 |
| 5 | 0.25 | 25 | 6.630 | -1.630 | 3.800 | 11.411 | 0.964 | 0.228 | 0.022 | 0.062 | 0.242 | 0.902 | -0.143 |
| | | 50 | 5.693 | -0.693 | 1.767 | 6.141 | 0.962 | 0.237 | 0.013 | 0.044 | 0.174 | 0.932 | -0.055 |
| | | 100 | 5.312 | -0.313 | 1.020 | 3.791 | 0.955 | 0.244 | 0.006 | 0.031 | 0.123 | 0.934 | -0.024 |
| 5 | 0.50 | 25 | 6.019 | -1.019 | 2.382 | 7.834 | 0.972 | 0.480 | 0.020 | 0.058 | 0.226 | 0.930 | -0.091 |
| | | 50 | 5.447 | -0.447 | 1.264 | 4.665 | 0.976 | 0.489 | 0.011 | 0.041 | 0.159 | 0.946 | -0.036 |
| | | 100 | 5.197 | -0.197 | 0.799 | 3.049 | 0.948 | 0.495 | 0.005 | 0.028 | 0.111 | 0.946 | -0.016 |
| 5 | 0.75 | 25 | 5.878 | -0.878 | 2.115 | 7.055 | 0.976 | 0.738 | 0.012 | 0.036 | 0.139 | 0.940 | -0.050 |
| | | 50 | 5.388 | -0.388 | 1.161 | 4.319 | 0.963 | 0.743 | 0.007 | 0.025 | 0.096 | 0.957 | -0.020 |
| | | 100 | 5.173 | -0.173 | 0.747 | 2.860 | 0.951 | 0.747 | 0.003 | 0.017 | 0.067 | 0.945 | -0.009 |

## 6. Empirical Data Analysis

In this section three data sets from the book of Reiss and Thomas [18] have been considered for modeling with $\text{DGUD}(\alpha, p)$. The first dataset 1 provides the annual maximum flood discharges of Feather River in Oroville, California, from 1902 to 1960 (Reiss and Thomas[18]), the second and the third data sets are respectively Changery's (1982) annual maximum wind speeds of tropical and non – tropical storms for Jacksonville, Florida (Reiss and Thomas [18]).

Absolute differences between the empirical (observed) distribution function and fitted (theoretical) distribution function at the sample points have been plotted in Figures 9 to highlight the closeness with which the $\text{DGUD}(\alpha, p)$ has been able to model the data.

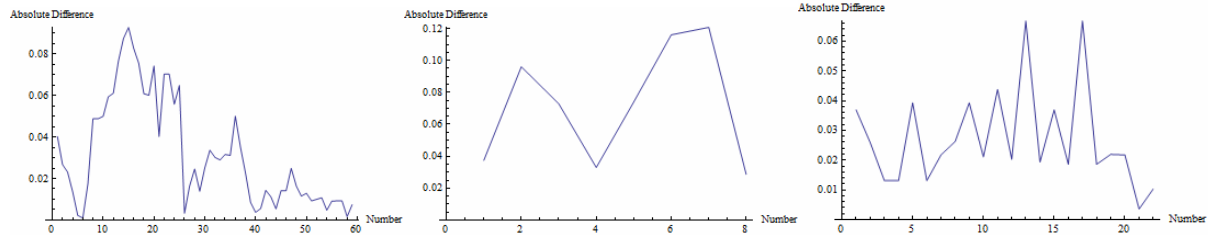

**Figure 9.** Absolute difference between the empirical cdf and theoretical cdf for Flood Data, Tropical Wind Data, Non-Tropical Wind Data

The Kolmogorov-Smirnov (KS) test for goodness of fit of discrete data [27, 28] has been used to check the performance of the model. The numerical findings are presented below in table 3.

**Table 3. Results of data fitting**

| Experiments | I. Annual Maximum Flood Discharges | II. Annual Maximum Wind Speeds (Tropical) | III. Annual Maximum Wind Speeds (Non-Tropical) |
|---|---|---|---|
| MLE | $\hat{\alpha} = 3.55393$ $\hat{p} = 0.99997$ | $\hat{\alpha} = 146.04395$ $\hat{p} = 0.89429$ | $\hat{\alpha} = 76275791$ $\hat{p} = 0.86040$ |
| Log L | -716.39427 | -30.57864 | -76.46989 |
| KS–test statistic value | 0.09274 | 0.12072 | 0.06679 |
| Lower bound of $p$ - value of | 0.56605 | 0.89961 | 0.60218 |

From the table 3, it can be seen that for all the data sets $DGUD(\alpha, p)$ has given good fit in terms of the discrete KS test statistics and corresponding *p*-values. Further the various cdf based statistics will return the same value as in case of adopting the continuous Gumbel since this discrete Gumbel retains the same cdf as that of continuous one.

## 7. Concluding Remarks

A new two parameters discrete analogue of Gumbel distribution with support on $Z$ has been proposed, it's important distributional, monotonic and reliability aspects have been investigated. A model adequacy check has been developed for the proposed distribution. Different estimation methods including the maximum likelihood estimation have been presented. Simulation study has been carried out to establish the accuracy and precision of the MLEs. Application of the proposed distribution in modeling three real life discrete extreme value data sets related to flood analysis and annual maximum wind speeds has been illustrated. It is envisaged that the proposed distribution defined $Z$ possessing many flexible properties will be a useful contribution to the field of count data modeling

**References**

[1] B. Mark, T. Bergstrom, Log-concave Probability and its Applications, Economic Theory 26 (2005) 445-469.

[2] C.D. Lai, Issues Concerning Constructions of Discrete Lifetime Models, Qual. Technol. Quant. Manag. 10, 2 (2013) 251-262.

[3] D. Roy, The Discrete Normal Distribution, Comm. Statist. Theory Methods, 32 (2003) 1871-1883.

[4] D. Roy, Discrete Rayleigh Distribution, IEEE Trans. Reliab. 53 (2004) 255-260.

[5] D. Roy, P.L. Gupta, Classifications of Discrete Lives, Microelectron. Reliab., 32, (1992) 1459-1473.

[6] F.W. Steutel, Log-concave and Log-convex Distributions, Encyclopedia of Statistical Sciences (eds S. Kotz, N.L. Johnson and C. B. Read), 5, 116-117. New York: Wiley, 1985.

[7] H. Krishna, P. S. Pundir, Discrete Maxwell Distribution, Interstat, http: // interstat. statjournals.net / YEAR / 2007 / articles / 0711003.pdf.

[8] H. Krishna, P.S. Pundir, Discrete Burr and Discrete Pareto distributions, Stat. Methodol. 6 (2009) 177-188.

[9] http://www.math.uah.edu/stat/special/ExtremeValue.html, The Extreme Value Distribution, accessed on 20.03.2014.

[10] J. Keilson, H. Gerber, Some Results for Discrete Unimodality, J. Amer. Statist. Assoc. 66, 334 (1971) 386-389.

[11] M. Xie, O. Gaudoin, C. Bracwuemond, Redefining Failure Rate Function for Discrete Distribution, J. Reliab. Qual. Safety Eng. 9 (2002) 275-285.

[12] M.A. Jazi, C.D. Lai, M.H. Alamatsaz, A Discrete Inverse Weibull Distribution and Estimation of its Parameters, Stat. Methodol. 7 (2010) 121-132.

[13] M.S.A. Khan, A. Khalique, A.M. Aboummoh, On Estimating Parameters in a Discrete Weibull Distribution, IEEE Trans. Reliab. 38 (1989) 348-350.

[14] N.L. Johnson, A.W. Kemp, S. Kotz, Continuous Univariate Distributions, Vol. 2, 2nd edition, Wiley, New York, 2004.

[15] N.L. Johnson, A.W. Kemp, S. Kotz, Univariate Discrete Distributions, 2nd edition, Wiley, New York, 2005.

[16] P.L. Gupta, R.C. Gupta, R.C. Tripathi, On the Monotonic Properties of Discrete Failure Rates, J. Statist. Plann. Inference 65 (1997) 255–268.

[17] R.C. Gupta, S.H. Ong, Analysis of Long-tailed Count Data by Poisson Mixtures, Comm. Statist. Theory Methods, 34 (2005) 557-573.

[18] R.D. Reiss, M. Thomas, Statistical Analysis of Extreme Events, 3rd ed., Birkhauser, Boston, Berlin, Germany, 2007.

[19] S. Chakraborty, A New Discrete Distribution Related to Generalized Gamma Distribution and its properties, Comm. Statist. Theory Methods, published on line, August, 2013.

[20] S. Chakraborty, D. Chakravarty, Discrete Gamma Distribution: Properties and Parameter Estimation, Comm. Statist. Theory Methods, 41(2012) 3301-3324.

[21] S. Chakraborty, D. Chakravarty, A New Discrete Probability Distribution with Integer Support on (-∞, ∞), Comm. Statist. Theory Methods, published on line, July, 2013.

[22] S. Das, P.K. Mandal, D. Ghosh, On Homogeneous Skeweness of Unimodal Distributions, Sankhya B 71 (2009) 187-205.

[23] S. Kotz, S. Nadarajah, Extreme Value Distributions: Theory and Applications, Imperial College Press, London, 2000.

[24] T. Nakagawa, S. Osaki, The Discrete Weibull Distribution, IEEE Trans. Reliab. 24 (1975) 300-301.

[25] V.K. Rohatgi, E.A.K. Md. Saleh, An Introduction to Probability and Statistics, 2nd edition, John-Wiley & Sons, New York, 2001.

[26] W.E. Stein, R. Dattero, A New Discrete Weibull Distribution, IEEE Trans. Reliab. 33 (1984) 196-197.

[27] W.J. Conover, A Kolmogorov Goodness-of-Fit Test for Discontinuous Distributions, J. Amer. Statist. Assoc. 67 339 (1972) 591-596.

[28] W.J. Conover, Practical Nonparametric Statistics, 3rd Edition, Wiley India, New Delhi, 1999.

[29] Y.A. Mark, Log-concave Probability Distributions: Theory and Statistical Testing, working paper, 96 - 01, Published by Center for Labour Market and Social Research, University of Aarhus and the Aarhus School of Business, Denmark, 1996.